# A New Method to Prove Goldbach Conjecture, Twin Primes Conjecture and Other Two Propositions

**Kaida Shi**

Department of Mathematics, Zhejiang Ocean University,
Zhoushan City 316004, Zhejiang Province, China

**Abstract**   By creating an identical method, the well-known world's baffling problems---Goldbach conjecture, twin primes conjecture and other two propositions have been proved.

**MR(2000): 11M**

In 1742, German mathematician C. Goldbach (1690~1764) in his letter addressed to Swiss mathematician L. Euler (1707~1783) wrote:

**Proposition (A)** Every even integer ($\geq 6$) is the sum of two odd primes;

**Proposition (B)** Every odd integer ($\geq 9$) is the sum of three odd primes.

They were called "Goldbach conjecture".

In 1900, German mathematician D. Hilbert (1862~1943) suggested 23 baffling mathematical problems to world's mathematical field at the 2nd International Congress of Mathematicians (held at Paris, France). The Goldbach conjecture and twin primes conjecture are a part of the 8-th problem[1,3].

In about 1920, English mathematicians G. H. Hardy, J. E. Littlewood and Indian mathematician S. Ramanujan raised the "Circle method"; Norway mathematician V. Brun raised the "Sieve method"; In 1933, Russian mathematician Шнирельман raised the "Density". These are all very important attempts along different directions for investigating Goldbach conjecture. Chinese mathematician Jingrun Chen improved Norway mathematician Selberg's "Sieve method", and obtained an important result {1, 2}[2]. All these carried powerfully the Analysis number theory and other many mathematical branches to a new development stage[1,3].

The aim of this paper is to create a new method, in order to prove Goldbach conjecture, twin primes conjecture and other two propositions raised by German mathematician E. Landua (1877~1938) on 5th International Congress of Mathematician held at Cambridge, England on 1912 [1,3].

## 1 The equivalent theorem of Goldbach—Vinogradov theorem and its proof

As is known to all, from the proposition (A), we can induce easily the Proposition (B). If the large even integer be denoted as $N$, so we have

$$N = p_1 + p_2.$$

If add 3 to two sides of above expression, we have

$$N + 3 = p_1 + p_2 + 3.$$

Because $N$ can take all large even integers, therefore, $N + 3$ will take all large odd integers. At the same time, because 3 is an odd prime, therefore, we can rewrite it as $p_3$. Thus, we have

$$N + 3 = p_1 + p_2 + p_3.$$

So. The Proposition (B) is proved. But, however, from the Proposition (B), we cannot prove the Proposition (A). This is because of when we represent a large odd integer as the sum of three odd primes, we cannot firm a certain prime from three primes.

For proving the Proposition (B), many mathematicians and amateurs have paid gigantic price. Fortunately, up to 1973, Chinese mathematician Jingrun Chen only obtained the result {1, 2} on the Proposition (A) [2]. This is a most distinguished achievement on this topic in Chinese, also in international mathematical field. But, it's very regret, there are some distinguished mathematicians have considered that for proving the Proposition (A), if we continue to use Chen's method (Sieve method), then we cannot obtain final success [3]. Hence, we must create a new method to finish the unfinished probe.

For proving the Proposition (A), I have undergone a long process of tempering. Now, I suggest the following:

**Definition** Suppose that there exist a multi-valued function $P(M) = e^{2\pi i f(M)}$ whose every value corresponds to an odd prime, where

$$M = \sum_{j=1}^{|\alpha|} p_j,$$

$|\alpha|$ is the number of the primes of a large integer $M$ is represented as the sum of odd primes. We can discover that all values of the multi-valued function $P(M)$ which corresponds to $p_j$ ($j = 1, 2, \cdots, |\alpha|$) will distribute evenly in the unit-circumference whose center is at the origin $O$ of the complex plane. Although $\alpha$ can be taken a positive integer or a negative integer, but $j$ must be taken a positive integer, because of it is represented as an ordinal number.

**Theorem 1 (The equivalent theorem of Goldbach --Vinogradov theorem)** *Let $N$ be a large odd integer, $P(N)$ be the number of primes of $N$ is represented as the sum of odd primes, then we have*

$$P(N) = e^{2\pi i(\frac{j}{3} + (N-3))} = e(\frac{j}{3} + (N-3))$$
$$= \cos(2\pi(\frac{j}{3} + (N-3))) + i\sin(2\pi(\frac{j}{3} + (N-3))) \qquad (*)$$

$$(j = 1, 2, 3)$$

*Proof* According to Goldbach—Vinogradov theorem, we have known that there exist an absolute constant $c_1$, it enables for every large odd integer $N$ can be represented as the sum of three odd primes. But, this $c_1$ is a very large number (about equals to $e^{e^{16.038}} > 10^{4,000,000}$).

Now, let's use the mathematical induction to prove the correctness of this **equivalent theorem**. So, the gap of the number less than $e^{e^{16.038}}$ will be filled.

When $N = 9$, we have

$$P(9) = e^{2\pi i(\frac{j}{3}+(9-3))} = e(\frac{j}{3}+(9-3))$$

$$= \cos(2\pi(\frac{j}{3}+6)) + i\sin(2\pi(\frac{j}{3}+6))$$

$$= \begin{cases} -\frac{1}{2}+\frac{\sqrt{3}}{2}i; \\ -\frac{1}{2}-\frac{\sqrt{3}}{2}i; = 3. \\ 1. \end{cases}$$

Suppose that when $N = m$, the expression (*) is also correct, namely,

$$P(m) = e^{2\pi i(\frac{j}{3}+(m-3))} = e(\frac{j}{3}+(m-3))$$

$$= \cos(2\pi(\frac{j}{3}+(m-3))) + i\sin(2\pi(\frac{j}{3}+(m-3)))$$

$$= \begin{cases} -\frac{1}{2}+\frac{\sqrt{3}}{2}i; \\ -\frac{1}{2}-\frac{\sqrt{3}}{2}i; = 3. \\ 1. \end{cases}$$

then when $N = m+2$, we have

$$P(m+2) = e^{2\pi i(\frac{j}{3}+(m+2-3))} = e^{2\pi i(\frac{j}{3}+(m-3))} \cdot e^{4\pi i}$$

$$= e(\frac{j}{3}+(m-3)) \cdot (\cos 4\pi + i\sin 4\pi)$$

$$= [\cos(2\pi(\frac{j}{3}+(m-3))) + i\sin(2\pi(\frac{j}{3}+(m-3)))] \cdot 1$$

$$= \begin{cases} -\frac{1}{2}+\frac{\sqrt{3}}{2}i; \\ -\frac{1}{2}-\frac{\sqrt{3}}{2}i; = 3. \\ 1. \end{cases}$$

From above, we can discover that the theorem we concerned only involved the unique variable $N$ (large odd integer), and doesn't involved any others.

    This proved that all large odd integers can be represented as the sum of three odd primes. Also, it explained that the result of this theorem is in keeping with what Goldbach—Vinogradov theorem requires of it. So, this theorem has been proved.

**Corollary** *For general large integer* $M(=N+\beta)$, *we can rewrite the expression* (*) *as:*

$$P(M) = e^{2\pi i(\frac{j}{3-\beta}+(M-3))} = e(\frac{j}{3-\beta}+(M-3))$$
$$= \cos(2\pi(\frac{j}{3-\beta}+(M-3))) + i\sin(2\pi(\frac{j}{3-\beta}+(M-3)))$$
$$= 3-\beta = \alpha,$$
$$(j = 1, 2, 3, \cdots, |\alpha|).$$

*where* $|\alpha|$ *is the number of the primes of a large integer* $M$ *is represented as the sum of odd primes.*

## 2. The proof of Goldbach conjecture, twin primes conjecture and other two propositions

### 2.1 The proof of Goldbach conjecture

**Theorem 2 (Goldbach conjecture)** *A large even integer can be represented as the sum of two odd primes.*

***Proof*** Because $N$ is a large odd integer, therefore, $N+1$ is a large even integer. Substituting $N+1$ into the expression (*), we obtain:

$$P(N+1) = e^{2\pi i(\frac{j}{3-1}+(N+1-3))} = e(\frac{j}{2}+(N-2))$$
$$= \cos(2\pi(\frac{j}{2}+(N-2))) + i\sin(2\pi(\frac{j}{2}+(N-2)))$$
$$= \begin{cases} -1; \\ 1. \end{cases} = 2.$$

This proved that a large even integer $N+1$ can be represented as the sum of two odd primes. So, the even integer Goldbach conjecture (Proposition (A), namely the {1, 1}) has been proved.

### 2.2 The proof of twin primes conjecture

**Theorem 3 (Twin primes conjecture)** *For every prime* $p$, *we can write an integer* $p+2$. *If* $p$ *takes all the odd primes, then we can obtain infinite integers whose form as* $p+2$, *and there exist infinite primes in these integers.*

***Proof*** Because $N$ is a large odd integer, then $N+2$ is also a large odd integer. Substituting it into the expression (*), we obtain

$$P(N+2) = e^{2\pi i(\frac{j}{3-2}+(N+2-3))} = e(\frac{j}{1}+(N-1))$$
$$= \cos(2\pi(\frac{j}{1}+(N-1))) + i\sin(2\pi(\frac{j}{1}+(N-1)))$$
$$= \{1. = 1.$$

In this moment, a large odd integer $N+2$ is only represented by an odd prime, namely, all the composite integer within all large odd integers $N+2$ have been sieved and remained integers are only odd primes. Because $N$ are infinite, therefore, these odd primes are also infinite, which more than and equal to 11.

Since $N+2$ are all odd primes $p$ which more than and equal to 11, then $N+4$ must equal to $p+2$, therefore, we have

$$P(N+4) = e^{2\pi i(\frac{j}{3-4}+(N+4-3))} = e(\frac{j}{-1}+(N-(-1)))$$
$$= \cos(2\pi(\frac{j}{-1}+(N-(-1)))) + i\sin(2\pi(\frac{j}{-1}+(N-(-1))))$$
$$= \{1. = 1.$$

In this moment, a large odd integer $N+4$ is also represented by an odd prime, namely, all composite integers within all large odd integers $N+4$ have been sieved and remained integers are only odd primes. Because $N$ are infinite, therefore, these odd primes are also infinite which more than and equal to 13. So, the twin primes conjecture has been proved.

**2. 3 The proofs of other two propositions**

**Theorem 4 (Proposition (C))** *There exist a positive integer $k$, it enables all integers $(\geq 2)$ can be represented as the sum of $k$ odd primes.*

***Proof*** We can verify concretely every integer which less than a large integer. Now, we consider only the situation of $N$ is a large integer. If we substitute $N-1$ (even integer), $N-2$ (odd integer), $N-3$ (even integer), ······ into the expression (*), we can obtain respectively:

$$P(N-1) = e^{2\pi i(\frac{j}{3-(-1)}+(N-1-3))} = e(\frac{j}{4}+(N-4))$$

$$= \cos(2\pi(\frac{j}{4}+(N-4))) + i\sin(2\pi(\frac{j}{4}+(N-4)))$$

$$= \begin{cases} i; \\ -1; \\ -i; \\ 1. \end{cases} = 4.$$

$$P(N-2) = e^{2\pi i(\frac{j}{3-(-2)}+(N-2-3))} = e(\frac{j}{5}+(N-5))$$

$$= \cos(2\pi(\frac{j}{5}+(N-5))) + i\sin(2\pi(\frac{j}{5}+(N-5)))$$

$$= \begin{cases} \frac{\sqrt{5}-1}{4} + \frac{\sqrt{10+2\sqrt{5}}}{4}i; \\ -\frac{\sqrt{5}-1}{4} + \frac{\sqrt{10-2\sqrt{5}}}{4}i; \\ -\frac{\sqrt{5}+1}{4} - \frac{\sqrt{10-2\sqrt{5}}}{4}i; = 5. \\ \frac{\sqrt{5}-1}{4} - \frac{\sqrt{10+2\sqrt{5}}}{4}i; \\ 1. \end{cases}$$

$$P(N-3) = e^{2\pi i(\frac{j}{3-(-3)}+(N-3-3))} = e(\frac{j}{6}+(N-6))$$

$$= \cos(2\pi(\frac{j}{6}+(N-6))) + i\sin(2\pi(\frac{j}{6}+(N-6)))$$

$$= \begin{cases} \frac{1}{2} + \frac{\sqrt{3}}{2}i; \\ -\frac{1}{2} + \frac{\sqrt{3}}{2}i; \\ -1; \\ -\frac{1}{2} - \frac{\sqrt{3}}{2}i; = 6. \\ \frac{1}{2} - \frac{\sqrt{3}}{2}i; \\ 1. \end{cases}$$

………………………………………….

This proved that a large odd integer can be represented respectively as the sums of three, five, seven, ⋯⋯ odd primes, and a large even integer can be represented respectively as the sums of two, four, six, ······ odd primes.

**Theorem 5 (Proposition 1)** *The primes whose form as $n^2 + 1$ are infinite.*

***Proof*** Because $n^2 + 1$ is an odd integer, therefore, $n^2 - 1$ is also an odd integer. For the odd integer $n^2 - 1$, we have

$$P(n^2 - 1) = e^{2\pi i(\frac{j}{3} + ((n^2-1)-3))} = e(\frac{j}{3} + ((n^2-1) - 3))$$

$$= \cos(2\pi(\frac{j}{3} + ((n^2-1) - 3))) + i\sin(2\pi(\frac{j}{3} + ((n^2-1) - 3)))$$

$$= \begin{cases} -\frac{1}{2} + \frac{\sqrt{3}}{2}i; \\ -\frac{1}{2} - \frac{\sqrt{3}}{2}i; = 3. \\ 1. \end{cases}$$

That is to say the odd integer $n^2 - 1$ can be represented as the sum of three odd primes, then for $n^2 + 1$, we have

$$P(n^2 + 1) = e^{2\pi i(\frac{j}{3-2} + ((n^2-1)+2-3))} = e(\frac{j}{1} + ((n^2-1) + 2 - 3))$$

$$= \cos(2\pi(\frac{j}{1} + (n^2 - 2))) + i\sin(2\pi(\frac{j}{1} + (n^2 - 2)))$$

$$= \{1. = 1.$$

In this moment, an odd integer $n^2 + 1$ is only represented by an odd prime, namely, all the composite integers within all odd integers $n^2 + 1$ have been sieved and remained integers are only the odd primes. Because $n$ are infinite, therefore, the odd primes whose form as $n^2 + 1$ are also infinite.

Thus, above two propositions raised by E. Landau have been proved.

## 3 Conclusion

The times since Goldbach conjecture was raised to now, it has undergone for 259 years. For proving this "unfathomable riddle", many mathematicians and amateurs of many countries have expended much energy. But, up to now, nobody can find a suitable method to carry on this research, enables for proving this problem is only depend on available theories. Obviously, the available theories are not very perfect, they cannot be adopted to prove the very complicated prime distribution problem, therefore they also cannot be used to prove Goldbach conjecture[1,3,4,5]. Hence, we must open actively a new road, find a new law, create a new method, in order to carry on new probe. So, to obtain success will be possible.

ZHEJIANG OCEAN INSTITUTE, ZJ, CHINA